\documentstyle[12pt]{article}
\newcommand{\be}{\begin{equation}}
\newcommand{\ee}{\end{equation}}
\newcommand{\bea}{\begin{eqnarray}}
\newcommand{\eea}{\end{eqnarray}}
\newcommand{\ba}{\begin{array}}
\newcommand{\ea}{\end{array}}

\newcommand{\Z}{Z\!\!\!\!Z}

\newcommand{\bc}{\begin{center}}
\newcommand{\ec}{\end{center}}
\newcommand{\ben}{\begin{enumerate}}
\newcommand{\een}{\end{enumerate}}
\newcommand{\bfi}{\begin{figure}}
\newcommand{\efi}{\end{figure}}

\newcommand{\bq}{\begin{quote}}
\newcommand{\eq}{\end{quote}}
\newcommand{\bqu}{\begin{quotation}}
\newcommand{\equ}{\end{quotation}}
\newenvironment{emphit}{\begin{itemize}}{\end{itemize}}
\newcommand{\bemp}{\begin{emphit}}
\newcommand{\eemp}{\end{emphit}}

\newcommand{\bt}{\begin{tabular}}
\newcommand{\et}{\end{tabular}}

\newtheorem{myth}{Theorem}[section]
\newtheorem{mylem}{Lemma}[section]
\newtheorem{mycor}{Corollary}[section]

\newtheorem{mydef}{Definition}[section]
\newtheorem{myrem}{Remark}[section]

\begin{document}
\date{}
\title{Algebraic properties of some varieties of central loops
}
\author{T\`em\'it\d {\'o}p\d {\'e} Gb\d {\'o}l\'ah\`an Jaiy\'e\d
ol\'a$^1$\\
 \&\\
John Ol\'us\d ol\'a Ad\'en\'iran$^2$\thanks{Corresponding author}}
\maketitle
\begin{abstract}
Isotopes of C-loops with unique non-identity squares are shown to be
both C-loops and A-loops. The relationship between C-loops and
Steiner loops is further studied. Central loops with the weak and
cross inverse properties are also investigated. C-loops are found to
be Osborn loops if every element in them are squares.
\end{abstract}

\section{Introduction}

C-loops are one of the least studied loops. Few publications that
have considered C-loops include Fenyves \cite{phd50}, \cite{phd56},
Beg \cite{phd169}, \cite{phd170}, Phillips et. al. \cite{phd9},
\cite{phd58}, \cite{phd59}, \cite{phd22}, Chein \cite{phd54} and
Solarin et. al. \cite{phd55}, \cite{phd52}, \cite{phd53},
\cite{phd10}. The difficulty in studying them is as a result of the
nature of their identities when compared with other Bol-Moufang
identities(the element occurring twice on both sides has no other
element separating it from itself). Latest publications on the study
of C-loops which has attracted fresh interest on the structure
include \cite{phd9}, \cite{phd58}, and \cite{phd59}.

LC-loops, RC-loops and C-loops are loops that satisfies the
identities  $(xx)(yz)=(x(xy))z$ , (zy)$(xx)=z((yx)x)$ and
$x(y(yz))=((xy)y)z$ respectively. Fenyves' work in \cite{phd56} was
completed in \cite{phd9}. Fenyves proved that LC-loops and RC-loops
are defined by three equivalent identities. But in \cite{phd9} and
\cite{phd61}, it was shown that LC-loops and RC-loops are defined by
four equivalent identities. Solarin \cite{phd53} named the fourth
identities left middle(LM-) and right middle(RM-) identities and
loops that obey them are called LM-loops and RM-loops respectively.
These terminologies were also used in \cite{phd51}. A C-loop is both
an LC-loop and an RC-loop (\cite{phd56}). Their basic properties are
found in \cite{phd58}, \cite{phd56} and \cite{phd49}.

The right and left translation maps on a loop $(L,\cdot )$ are
bijections
\begin{displaymath}
R_x~:~L\to L~\textrm{and}~L_x~:~L\to L~\textrm{defined
as}~yR_x=y\cdot x~\textrm{and}~yL_x=x\cdot y~\textrm{respectively
for all} ~x,y\in L.
\end{displaymath}

\begin{mydef}
Let $(L, \cdot )$ be a loop. The left nucleus of $L$ is the set
\begin{displaymath}
N_\lambda (L, \cdot )=\{a\in L : ax\cdot y=a\cdot xy~\forall~x, y\in
L\}.
\end{displaymath}
The right nucleus of $L$ is the set
\begin{displaymath}
N_\rho (L, \cdot )=\{a\in L : y\cdot xa=yx\cdot a~\forall~ x, y\in
L\}.
\end{displaymath}
The middle nucleus of $L$ is the set
\begin{displaymath}
N_\mu (L, \cdot )=\{a\in L : ya\cdot x=y\cdot ax~\forall~x, y\in
L\}.
\end{displaymath}
The nucleus of $L$ is the set
\begin{displaymath}
N(L, \cdot )=N_\lambda (L, \cdot )\cap N_\rho (L, \cdot )\cap N_\mu
(L, \cdot ).
\end{displaymath}
The centrum of $L$ is the set
\begin{displaymath}
C(L, \cdot )=\{a\in L : ax=xa~\forall~x\in L\}.
\end{displaymath}
The center of $L$ is the set
\begin{displaymath}
Z(L, \cdot )=N(L, \cdot )\cap C(L, \cdot ).
\end{displaymath}
\end{mydef}
$L$ is said to be a centrum square loop if~ $x^2\in C(L, \cdot )$
for all $x\in L$. $L$ is said to be a central square loop if~
$x^2\in Z(L, \cdot )$ for all $x\in L$. $L$ is said to be left
alternative if for all $x, y\in L,~ x\cdot xy=x^2y$ and is said to
right alternative if for all $x, y\in L,~ yx\cdot x=yx^2$. Thus, $L$
is said to be alternative if it is both left and right alternative.
The triple $(U, V, W)$ such that $U, V, W\in SYM(L, \cdot )$ is
called an autotopism of $L$ if and only if
\begin{displaymath}
xU\cdot yV=(x\cdot y)W~\forall ~x, y\in L.
\end{displaymath}
$SYM(L, \cdot )$ is called the permutation group of the loop
$(L,\cdot )$. The group of autotopisms of $L$ is denoted by $AUT(L,
\cdot )$. Let $(L, \cdot )$ and $(G, \circ )$ be two distinct loops.
The triple $(U, V, W) : (L, \cdot )\to (G, \circ )$ such that $U, V,
W : L\to G$ are bijections is called a loop isotopism if and only if
\begin{displaymath}
xU\circ yV=(x\cdot y)W~\forall ~x, y\in L.
\end{displaymath}

We investigate central loops with the unique non-identity
commutators, associators and squares and those with unique
non-trivial squares are also found to be C-loops whose isotopes are
C-loops and A-loops. The relationship between C-loops and Steiner
loops is further studied.

Central loops with the weak and cross inverse properties are also
investigated. C-loops are found to be Osborn loops if every
element in them are squares, hence the possibility for a
non-associative Osborn C-loop to exist is envisaged.

For definition of concepts in theory of loops readers may consult
\cite{phd41}, \cite{phd51} and \cite{phd3}.

\section{Preliminaries}

\begin{mydef}(\cite{phd42})
Let $a,b$ and $c$ be three elements of a loop $L$. The loop
commutator of $a$ and $b$ is the unique element $(a,b)$ of $L$
which satisfies $ab=(ba)(a,b)$ and the loop associator of $a,b$
and $c$ is the unique element $(a,b,c)$ of $L$ which satisfies
$(ab)c=\{a(bc)\}(a,b,c)$.
\end{mydef}

\paragraph{}
If $X,Y,$ and $Z$ are subsets of a loop $L$, we denote by $(X,Y)$
and $(X,Y,Z)$, respectively, the set of all commutators of the
form $(x,y)$ and all the associators of the form $(x,y,z)$, where
$x\in X,y\in Y,z\in Z$.

\begin{mydef}(\cite{phd42})
A 'unique non-identity commutator' is an element $s\ne e$($e$ is the
identity element) in a loop $L$ with the property that every
commutator in $L$ is $e$ or $s$.

A 'unique non-identity commutator associator' is an element $s\ne
e$ in a loop $L$ with the property that every commutator in $L$ is
$e$ or $s$ and every associator is $e$ or $s$.

A 'unique non-identity square' or 'non-trivial square' is an element
$s\ne e$ in a loop $L$ with the property that every square in $L$ is
$e$ or $s$.
\end{mydef}

\begin{mydef}\label{0:0.6}
A loop$(L,\cdot)$ is called a weak inverse property loop (W. I. P.
L.) if and only if it obeys the weak inverse property (W. I. P.) :
$y(xy)^\rho=x^\rho$ for all $x,y\in L$ while $L$ is called a cross
inverse property loop (C. I. P. L.) if and only if it obeys the
cross inverse property (C. I. P.) : $xy\cdot x^\rho=y$. $(L,\cdot)$
is called a left inverse property loop or right inverse property
loop (L.I.P.L. or R.I.P.L.) if and only if it obeys the left inverse
property or right inverse property(L.I.P or R.I.P): $ x^\lambda (xy)
= y ~\textrm{or}~(yx) x^\rho=y$. Hence, it is called an inverse
property loop (I.P.L.) if and only if it has the inverse property
(I.P.) i.e. it is both an L.I.P. and an R.I.P. loop.
\end{mydef}

Most of our results and proofs, are stated and written in dual
form relative to RC-loops and LC-loops. That is, a statement like
'LC(RC)-loop... A(B)' where 'A' and 'B' are some equations or
expressions simply means 'A' is for LC-loops while 'B' is for
RC-loops. This is done so that results on LC-loops and RC-loops
can be combined to derive those on C-loops. For instance an
LC(RC)-loop is a L.I.P.L.(R.I.P.L) loop while a C-loop in an
I.P.L. loop.

\section{Inner mappings}
\begin{mylem}\label{c:auto2}
Let $L$ be a C-loop. Then for each $(A,B,C)\in AUT(L,\cdot)$, there
exists a unique pair of $(S_1,T_1,R_1), (S_2,T_2,R_2) \in
AUT(L,\cdot)$ for each $x\in L$ such that $L_x^2=S_2^{-1}S_1,
R_x^2=T_1^{-1}T_2, R_x^{-2}L_x^2=R_2^{-1}R_1,
R_1^{-1}R_2T_2^{-1}T_1S_2^{-1}S_1=I$.
\end{mylem}
{\bf Proof}\\If $L$ is a C-loop, then
$(L_x^2,I,L_x^2),(I,R_x^2,R_x^2)\in AUT(L)$ for all $x\in L$. So
there exist $(S_1,T_1,R_1), (S_2,T_2,R_2) \in AUT(L,\cdot)$ such
that
\begin{displaymath}
(S_1,T_1,R_1)=(A,B,C)(L_x^2,I,L_x^2)\in AUT(L)
\end{displaymath}
\begin{displaymath}
(S_2,T_2,R_2)=(A,B,C)(I,R_x^2,R_x^2)\in AUT(L).
\end{displaymath}
Hence, the conditions hold although the identities do not depend on
$(A,B,C)$, but the uniqueness does.$\spadesuit$

\begin{myth}\label{c:auto2comp}
Let $L$ be a C-loop and let there exist a unique pair of autotopisms
$(S_1,T_1,R_1), (S_2,T_2,R_2)$ such that the conditions
$L_x^2=S_2^{-1}S_1, R_x^2=T_1^{-1}T_2$ and
$R_x^{-2}L_x^2=R_2^{-1}R_1$ hold for each fixed $x\in L$. If
$\alpha_1=S_1^{-1}, \alpha_2=S_2^{-1},
\beta_1=T_1^{-1},\beta_2=T_2^{-1}, \gamma_1=R_1^{-1}$ and
$\gamma_2=R_2^{-1}$, then:
\begin{displaymath}
(x^2y)\alpha_1 =y\alpha_2,\qquad (yx^2)\beta_2 =y\beta_1,\qquad
(x^2yx^{-2})\gamma_1 =y\gamma_2~\forall~x,y\in L.
\end{displaymath}
\end{myth}
{\bf Proof}\\
From Lemma~\ref{c:auto2}:
\begin{displaymath}
L_x^2=S_2^{-1}S_1,R_x^2=T_1^{-1}T_2,R_x^{-2}L_x^2=R_2^{-1}R_1.
\end{displaymath}
Keeping in mind that a C-loop is power associative and nuclear
square, we have the following proofs.
\begin{enumerate}
\item $L_x^2=S_2^{-1}S_1$ implies $yL_x^2=yS_2^{-1}S_1$ for all $y\in L$ implies $yL_{x^2}=yS_2^{-1}S_1$ implies
$x^2y=yS_2^{-1}S_1$ implies $(x^2y)S_1^{-1}=yS_2^{-1}$ implies
$x^2y\alpha_1=y\alpha_2$.
\item $R_x^2=T_1^{-1}T_2$ implies $yR_x^2=yT_1^{-1}T_2$ for all $y\in L$ implies $yx^2=yT_1^{-1}T_2$ implies
$yx^2T_2^{-1}=yT_1^{-1}$ implies $yx^2\beta=y\beta_1$.
\item $R_x^{-2}L_x^2=R_2^{-1}R_1$ implies $yR_x^{-2}L_x^2=yR_2^{-1}R_1$ for all $y\in L$ implies
$x^2yx^{-2}=yR_2^{-1}R_1$ implies $(x^2yx^{-2})R_1^{-1}=yR_2^{-1}$
implies $(x^2yx^{-2})\gamma_1=y\gamma_2$.
\end{enumerate}
$\spadesuit$

\begin{mycor}\label{c:autotopism}
Let $L$ be a C-loop. An autotopism of $L$ can be constructed if
there exists at least an $x\in L$ such that $x^2\neq e$. The inverse
can also be constructed.
\end{mycor}
{\bf Proof}\\ We need Lemma~\ref{c:auto2} and
Theorem~\ref{c:auto2comp}. If $x^2=e$, then the autotopism is
trivial. Since $L$ is a C-loop, using Lemma~\ref{c:auto2} and
Theorem~\ref{c:auto2comp}, it will be noticed that $(\alpha_1
S_2,\beta_1 T_2,\gamma_1 R_2)\in AUT(L)$ and $(\alpha_2 S_1,\beta_2
T_1,\gamma_2 R_1)=(\alpha_1 S_2,\beta_1 T_2,\gamma_1 R_2)^{-1}$.
Hence the proof. $\spadesuit$

\begin{mylem}\label{c:g2r1}
Let $(L,\cdot )$ be a C-loop. Then the third component mapping
$\gamma_2R_1~:~L\to L$ defined by $y\gamma_2R_1=x^2yx^{-2}$ for all
$x\in L$ of the constructed autotopism $(\alpha_2 S_1,\beta_2
T_1,\gamma_2 R_1)\in AUT(L)$ is :
\begin{enumerate}
\item an automorphism, \item a semi-automorphism, \item a middle
inner mapping, \item a pseudo-automorphism with companion $x^2$.
\end{enumerate}
\end{mylem}
{\bf Proof} \\
\begin{enumerate}
\item $\gamma_2R_1$ is already a bijection by the construction of
the autotopism $(\alpha_2 S_1,\beta_2 T_1,\gamma_2 R_1)\in AUT(L)$.
So we need only to show that it is an homomorphism. Let $y_1,y_2\in
L$, then :
$(y_1y_2)\gamma_2R_1=(x^2y_1x^{-2})(x^2y_2x^{-2})=y_1\gamma_2R_1\cdot
y_2\gamma_2R_1$. Whence, $\gamma_2R_1$ is an automorphism.
\item It will be observed that, $e\gamma_1=e\gamma_2$, hence
$e\gamma_2R_1=e$. $(zy\cdot z)\gamma_2R_1=x^2(zy\cdot
z)x^{-2}=x^2((zy\cdot z)x^{-2}) =(x^2zx^{-2})(x^2yx^{-2})\cdot
z\gamma_2R_1 =(z\gamma_2R_1\cdot y\gamma_2R_1)\cdot z\gamma_2R_1$.
Thence, $\gamma_2R_1$ is a semi-automorphism
\item We have already settled the fact that $e\gamma_2R_1=e$.
$y\gamma_2R_1=yR_{x^{-2}}L_{(x^{-2})^{-1}}=yT(x^{-2})$ for all $y\in
L$ implies $\gamma_2R_1=T(x^{-2})\in Inn(L)$. Hence $\gamma_2R_1$ is
a middle inner mapping.
\item This is true from (1) by using the
fact that says an automorphism in a C-loop $L$ is a
pseudo-automorphism with companion $x^2$ for all $x\in L$.
\end{enumerate}
$\spadesuit$

\begin{mylem}\label{c:innermap}
Let $(L,\cdot)$ be a C-loop. Then :
\begin{enumerate}
\item
$T(x^{-1})=R_xT(x^{-2})L_x^{-1},T(x)^2=R_xT(x^{-1})^{-1}L_x^{-1}$,
\item
$T(x^n)=R_x^{n-1}T(x)L_x^{1-n},T(x^{-n})=R_x^{1-n}T(x^{-1})L_x^{n-1}$
for all $n\in {\Z}^+$, \item $R(x,x)=I,L(x,x)=I$.
\end{enumerate}
\end{mylem}
{\bf Proof} \\
\begin{enumerate}
\item Using the mapping $\gamma_2R_1$ in Lemma~\ref{c:g2r1},
defined by $y\gamma_2R_1=x^2yx^{-2}$ for all $x,y\in L$,
$y\gamma_2R_1=x^2yx^{-2}=yR_{x^{-2}}L_{x^2}=yR_x^{-1}R_x^{-1}L_xL_x
=yR_x^{-1}T(x^{-1})L_x$. Thus, $\gamma_2R_1=R_x^{-1}T(x^{-1})L_x$.
In Lemma~\ref{c:g2r1}, we showed that $\gamma_2R_1$ is the middle
inner mapping $T(x^{-2})$. Hence, $T(x^{-2})=R_x^{-1}T(x^{-1})L_x$
implies $T(x^{-1})=R_xT(x^{-2})L_x^{-1}$.
$T(x)^2=R_xL_x^{-1}R_xL_x^{-1}=R_x(R_{x^{-1}}L_{x^{-1}}^{-1})^{-1}L_x^{-1}
=R_xT(x^{-1})^{-1}L_x^{-1}$.
\item We prove these by induction. When
\begin{displaymath}
n=1,~
T(x)=R_x^{1-1}T(x)L_x^{1-1}=R_{x^0}T(x)L_{x^0}=T(x)~\textrm{for
all}~x\in L,
\end{displaymath}
\begin{displaymath}
n=2,~
T(x^2)=T(xx)=R_{x^2}L_{x^2}^{-1}=R_xR_xL_x^{-1}L_x^{-1}=R_xT(x)L_x^{-1}~\textrm{for
all}~x\in L,
\end{displaymath}
\begin{displaymath}
n=3,~
T(x^3)=T(x^2x)=R_{x^2x}L_{(x^2x)^{-1}}=R_{x^2}R_xL_{x^{-1}x^{-2}}
=R_{x^2}R_xL_{x^{-1}}L_{x^{-2}}
\end{displaymath}
\begin{displaymath}
=R_x^2T(x)L_x^{-2}~\textrm{for all}~x\in L.
\end{displaymath}

Let $n=k,~ T(x^k)=R_x^{k-1}T(x)L_x^{1-k}$.

If
$n=k+1,~T(x^{k+1})=T(x^{k-1}x^2)=R_{x^{k-1}x^2}L_{(x^{k-1}x^2)}^{-1}=R_{x^{k-1}x^2}L_{x^{-2}x^{1-k}}
=R_{x^{k-1}}R_{x^2}L_{x^{-2}}L_{x^{1-k}}=R_{x^{k-1}}T(x^2)L_{x^{1-k}}
=R_x^{k-1}R_xT(x)L_x^{-1}L_x^{1-k}=R_x^kT(x)L_x^{-k}$. Thus,
$T(x^n)=R_x^{n-1}T(x)L_x^{1-n}$ for all $n\in {\Z}^+$. Replace $x$
by $x^{-1}$, then
$T(x^{-n})=T(({x^{-1}})^n)=R_{x^{-1}}^{n-1}T(x^{-1})L_{x^{-1}}^{1-n}=R_x^{1-n}T(x^{-1})L_x^{n-1}$.
Thus, $T(x^{-n})=R_x^{1-n}T(x^{-1})L_x^{n-1}$ for all $n\in {\Z}^+$.
\item The right and left inner mappings as defined gives :
$R(x,x)=R_x^2R_x^{-2}=I,L(x,x)=L_x^2L_x^{-2}=I$.
\end{enumerate} $\spadesuit$

\begin{myrem}
Lemma~\ref{c:g2r1} gives an example of a bijective mapping that is
an automorphism, pseudo-automorphism, semi-automorphism and an
inner mapping.
\end{myrem}

\section{Relationship between C-loops and Steiner loops}

For a loop $(L,\cdot)$, the bijection $J~:~L\to L$ is defined by
$xJ=x^{-1}$ for all $x\in L$.
 $L$ is called a Steiner loop if and only if
\begin{displaymath}
x^2=e~, ~yx\cdot x=y~\textrm{and}~xy=yx~\forall ~x, y\in L.
\end{displaymath}

\begin{myth}\label{c:exponent4}
In a C-loop $(L,\cdot )$, if any of the following is true for all
$z\in L$:
\begin{enumerate}
\item $(I,L_z^2,JL_z^2J)\in AUT(L),$ \item $(R_z^2,I,JR_z^2J)\in
AUT(L),$
\end{enumerate}
then, $L$ is a loop of exponent $4$.
\end{myth}
{\bf Proof} \\
\begin{enumerate}
\item If $(I,L_z^2,JL_z^2J)\in AUT(L)$ for all $z\in L$, then : $x\cdot
yL_z^2=(xy)JL_z^2J$ for all $x,y,z\in L$ implies $x\cdot
z^2y=xy\cdot z^{-2}$ implies $z^2y\cdot z^2=y$. Then $y^4=e$. Hence
$L$ is a C-loop of exponent $4$.
\item If $(R_z^2,I,JR_z^2J)\in AUT(L)$ for all $z\in L$, then :
$xR_z^2\cdot y=(xy)JR_z^2J$ for all $x,y,z\in L$ implies
$(xz^2)\cdot y=[(xy)^{-1}z^2]^{-1}$ implies $(xz^2)\cdot
y=z^{-2}(xy)$ implies $(xz^2)\cdot y=z^{-2}x\cdot y$ implies
$xz^2=z^{-2}x$ implies $z^4=e$. Hence $L$ is a C-loop of exponent
$4$.
\end{enumerate}
$\spadesuit$

\begin{myth}\label{c:centralsquare}
In a C-loop $L$, if the following are true for all $z\in L$ :
\begin{enumerate}
\item $(I,L_z^2,JL_z^2J)\in AUT(L),$ \item $(R_z^2,I,JR_z^2J)\in
AUT(L),$
\end{enumerate}
then, $L$ is a central square C-loop of exponent 4.
\end{myth}
{\bf Proof} \\
By the first hypothesis, If $(I,L_z^2,JL_z^2J)\in AUT(L)$ for all
$z\in L$, then : $x\cdot yL_z^2=(xy)JL_z^2J$ for all $x,y,z\in L$
implies $x\cdot z^2y=xy\cdot z^{-2}$.

By the second hypothesis, If $(R_z^2,I,JR_z^2J)\in AUT(L)$ for all
$z\in L$,then : $xR_z^2\cdot y=(xy)JR_z^2J$ for all $x,y,z\in L$
implies $xz^2\cdot y=z^{-2}(xy)$.

Using the two results above and keeping in mind that $L$ is a C-loop
we have :

$x\cdot z^2y=xz^2\cdot y$ if and only if $xy\cdot z^{-2}=z^{-2}\cdot
xy$. Let $t=xy$ then $tz^{-2}=z^{-2}t$ if and only if
$z^2t^{-1}=t^{-1}z^2$. Let $s=t^{-1}$ then $z^2\in C(L,\cdot )$ for
all $z\in L$.

Since $s$ is arbitrary in $L$, then the last result shows that $L$
is centrum square. Furthermore, C-loops have been found to be
nuclear square in \cite{phd58}, thus $z^2\in Z(L,\cdot )$. Hence $L$
is a central square C-loop. Finally, by Theorem~\ref{c:exponent4},
$x^4=e$.$\spadesuit$

\begin{mycor}\label{1:9}
In a C-loop $(L,\cdot )$, if $(I,L_z^2,JL_z^2J)\in AUT(L),$ and
$(R_z^2,I,JR_z^2J)\in AUT(L)$ for all $z\in L$, then the following
are true :
\begin{enumerate}
\item $L$ is flexible. \item $(xy)^2=(yx)^2$ for all $x,y\in L$.
\item $x\mapsto x^3$ is an anti-automorphism.
\end{enumerate}
\end{mycor}
{\bf Proof}\\
This follows by Theorem~\ref{c:centralsquare}, Lemma~5.1 and
Corollary~5.2 of \cite{phd59}.$\spadesuit$

\begin{myth}\label{1:10}
A central square C-loop of exponent 4 is a group.
\end{myth}
{\bf Proof}\\
To prove this, it shall be shown that the right inner mapping

$R(x,y)=I$ for all $x,y\in L$. Corollary~\ref{1:9} is used. Let
$w\in L$.

$wR(x,y)=wR_xR_yR_{xy}^{-1}=(wx)y\cdot (xy)^{-1}=(wx)(x^2yx^2)\cdot
(xy)^{-1}=(wx^3)(yx^2)\cdot (xy)^{-1}=(w^2(w^3x^3))(yx^2)\cdot
(xy)^{-1}=(w^2(xw)^3)(yx^2)\cdot (xy)^{-1}=w^2(xw)^3\cdot
(yx^2)(xy)^{-1}=w^2(xw)^3\cdot [y\cdot x^2(xy)^{-1}]=w^2(xw)^3\cdot
[y\cdot x^2(y^{-1}x^{-1})]=w^2(xw)^3\cdot [y(y^{-1}x^{-1}\cdot
x^2)]=w^2(xw)^3\cdot [y(y^{-1}x)]=w^2(xw)^3\cdot x=w^2(w^3x^3)\cdot
x=w^2\cdot (w^3x^3)x=w^2\cdot (w^3x^{-1})x=w^2w^3=w^5=w$ if and only
if $R(x,y)=I$ if and only if $R_xR_yR_{xy}^{-1}=I$ if and only if
$R_xR_y=R_{xy}$ if and only if $zR_xR_y=zR_{xy}$ if and only if
$zx\cdot y=z\cdot xy$ if and only if $L$ is a group. Hence the claim
is true.$\spadesuit$

\begin{mycor}\label{1:11}
In a C-loop $(L,\cdot )$, if $(I,L_z^2,JL_z^2J)\in AUT(L),$ and
$(R_z^2,I,JR_z^2J)\in AUT(L)$ for all $z\in L$, then $L$ is a group.
\end{mycor}
{\bf Proof}\\
This follows from Theorem~\ref{c:centralsquare} and
Theorem~\ref{1:10}.$\spadesuit$

\begin{myrem}
Central square C-loops of exponent 4 are A-loops.
\end{myrem}

\begin{myth}\label{c:st1}
Let $L$ be a C-loop. $L$ is a central square loop if and only if
$\gamma_2R_1=I$.
\end{myth}
{\bf Proof}\\
Let $\gamma_2R_1=I$ if and only if $T(x^{-2})=I$ for all $x\in L$ if
and only if $R_{x^{-2}}L_{x^2}=I$ if and only if $yx^2=x^2y$ if and
only if $L$ is central square.$\spadesuit$

\begin{myth}\label{c:st2}
Let $L$ be a C-loop such that the mapping $x\mapsto T(x)$ is a
bijection, then $L$ is of exponent $2$ if and only if
$\gamma_2R_1=I$.
\end{myth}
{\bf Proof}\\
Let $\gamma_2R_1=I$ if and only if $T(x^{-2})=I$ for all $x\in L$ if
and only if $T(x^{-2})=I=R_x^{-1}T(x^{-1})L_x$ by
Lemma~\ref{c:innermap} if and only if $T(x^{-1})=T(x)$ if and only
if $x^{-1}=x$ since $x\mapsto T(x)$ is a bijection if and only if
$L$ is a loop of exponent 2. $\spadesuit$

\begin{mycor}\label{c:important1}
A C -loop in which $x\mapsto T(x)$ is a bijection is of exponent 2
if and only if it is central square.
\end{mycor}
{\bf Proof} \\
This follows by Theorem~\ref{c:st1} and Theorem~\ref{c:st2}.
$\spadesuit$

\begin{myrem}
This Corollary leads to answer one of the questions we wish to
answer in this work. This we give below.
\end{myrem}

\begin{mycor}\label{c:important2}
In a central square C-loop $L$, the map $x\mapsto T(x)$ is a
bijection implies $L$ is a Steiner loop.
\end{mycor}
{\bf Proof} \\
By the converse of Corollary~\ref{c:important1}, a C-loop in which
$x\mapsto T(x)$ is a bijection, is of exponent 2 if it is central
square. By the result in \cite{phd58}, that an inverse property
loop of exponent 2 is exactly a Steiner loop and the fact that
C-loops are inverse property loops(\cite{phd58}), $L$ is a Steiner
loop. $\spadesuit$

\begin{mycor}
In a C-loop $(L,\cdot)$, if
\begin{enumerate}
\item $x\mapsto T(x)$ is a bijection, \item $(I,L_z^2,JL_z^2J)\in
AUT(L)$ and \item $(R_z^2,I,JR_z^2J)\in AUT(L)$ for all $z\in L$,
\end{enumerate}
then $L$ is a Steiner loop of exponent 4.
\end{mycor}
{\bf Proof} \\
According to Theorem~\ref{c:centralsquare}, $L$ is a central
square loop. Hence, since by hypothesis, $x\mapsto T(x)$ is a
bijection, then by Corollary~\ref{c:important2}, $L$ is a Steiner
loop. By Theorem~\ref{c:exponent4}, it has a an exponent of 4.
$\spadesuit$

\begin{mycor}\label{c:important3}
Let $(L,\cdot)$ be a C-loop such that the mapping $x\mapsto T(x) $
is a bijection. $L$ is a Steiner loop if and only if $L$ is a
central square C-loop.
\end{mycor}
{\bf Proof} \\
It was shown in \cite{phd58} that Steiner loops are C-loops. Recall
that Steiner loops are loops of exponent two by definition, hence by
Corollary~\ref{c:important1}, $L$ is central square since in $L$,
the mapping $x\mapsto T(x)$ is a bijection. Conversely : By
Corollary~\ref{c:important1}, a central square C-loop $L$ in which
the mapping $x\mapsto T(x)$ is a bijection is of exponent two.
Whence using the fact in \cite{phd58} that an inverse property loop
of exponent two is a steiner loop, the proof is complete.
$\spadesuit$

\subsection{Flexibility in C-loops}

Condition(s) under which a C-loop is flexible is given below.

\begin{mylem}\label{flex:onto}
A C-loop is flexible if the mapping $x\mapsto x^2$ is onto.
\end{mylem}
{\bf Proof} \\
Let $L$ be a C-loop . Then $yx^2\cdot y=y\cdot x^2y$ for all $x,y\in
L$. Thus, $L$ is square flexible, hence by the observation in
\cite{phd63} , $L$ is flexible since the mapping $x\mapsto x^2$ is
onto.

\begin{myth}\label{c:flex}
A C-loop $L$ is flexible if $(I,L_z^2,JL_z^2J),(R_z^2,I,JR_z^2J) \in
AUT(L)$ for all $z\in L$ and the middle inner mappings are of order
$2$.
\end{myth}
{\bf Proof} \\
By Lemma~\ref{c:innermap},

$T(x)^2=R_xT(x^{-1})^{-1}L_x^{-1}=R_x(R_xT(x^{-2})L_x^{-1})^{-1}L_x^{-1}=R_x(L_x(R_xT(x^{-2}))^{-1})L_x^{-1}
=R_x(L_xT(x^{-2})^{-1}R_x^{-1})L_x^{-1}=R_xL_xT(x^{-2})^{-1}R_x^{-1}L_x^{-1}=R_xL_xT(x^{-2})^{-1}(L_xR_x)^{-1}$.
So $T(x)^2=R_xL_xT(x^{-2})^{-1}(L_xR_x)^{-1}$ if and only if
$T(x)^2L_xR_x=R_xL_xT(x^{-2})^{-1}=R_xL_x(\gamma_2R_1)^{-1}=R_xL_x\gamma_1R_2$
if and only if $T(x)^2L_xR_x=R_xL_x\gamma_1R_2$. If
$|T(x)|=2,T(x)^2=I$ and if $\gamma_1R_2=I$ if and only if $L$ is
central square then $L_xR_x=R_xL_x$ if and only if $xy\cdot x=x\cdot
yx$ is a flexible loop. $\spadesuit$

\paragraph{}
Philips and Vojt\v echovsk\'y \cite{phd58} studied the close
relationship between C-loops and Steiner loops. In \cite{phd3}, it
is shown that Steiner loops are exactly commutative inverse
property loops of exponent 2. But in \cite{phd58}, this fact was
improved, so that commutativity is not a sufficient condition for
an inverse property loop of exponent 2 to be a Steiner loop. So
they said `Steiner loops are exactly inverse property loops of
exponent 2'. This result is general for inverse property loops
among which are C-loops. They also proved that Steiner loops are
C-loops. Corollary~\ref{c:important2} give conditions(apart from
exponent of 2) under which a C-loop is a Steiner loop while
Corollary~\ref{c:important3} give necessary and sufficient
conditions for C-loops to be Steiner loops.

It has been found out in \cite{phd3} that flexibility is possible in
a C-loop if the loop is commutative or diassociative. But C-loops
naturally do not even satisfy the latter. Apart from the condition
stated in Lemma~\ref{flex:onto}, Theorem~\ref{c:flex} when compared
with Corollary~5.2 of \cite{phd59} shows that some middle
inner-mappings do not need to be of exponent of $2$ for a C-loop to
be flexible.

\section{Unique non-identity Commutator Associator Square}

\begin{mylem}\label{1:12}
If $s$ is a unique non-identity commutator in a C-loop $L$,
$|s|=2$, $s\in C(L)$ and $s\in Z(L^2)$.
\end{mylem}
{\bf Proof}\\
$xy=(yx)(x,y)\Leftrightarrow
(x,y)=(yx)^{-1}(xy)=(x^{-1}y^{-1})(xy)$.
$(x,y)^{-1}=[(x^{-1}y^{-1})(xy)]^{-1}=(xy)^{-1}(x^{-1}y^{-1})^{-1}=(y^{-1}x^{-1})(yx)=(y,x)$.
Thus, $s^{-1}=s$ or $s^{-1}=e$ implies $s^2=e$ or $s=e$ implies
$s^2=e$.

If $xs\ne sx$, then $xs=(sx)s$ implies $x=sx$ implies $s=e$ implies
$xs=sx$ implies $s\in C(L)$. Hence, $s\in Z(L^2)$. $\spadesuit$

\begin{mylem}\label{1:13}
If $s$ is a unique non-identity associator in a C-loop $L$, $s\in
N(L)$.
\end{mylem}
{\bf Proof}\\
If $(xy)s\ne x(ys)$, then $(xy)s=x(ys)\cdot s$ implies $xy=x\cdot
ys$ implies $y=ys$ implies $s=e$ implies $(xy)s=x(ys)$ implies $s\in
N(L)$. $\spadesuit$

\begin{mylem}\label{c:unic}
Let $(L,\cdot)$ be a C-loop which has a unique non-identity
commutator associator $s$. Then $s$ is a central element of order 2.
\end{mylem}
{\bf Proof} \\
We shall keep in mind that $L$ as a C-loop has the inverse property.
$s\in (L,L)$ implies $s^{-1}\in (L,L)$ implies $s^{-1}=s$ since
$s^{-1}\ne e$ if and only if $s\ne e$, hence $s^2=e$. Let $xs\ne sx$
for some $x,y\in L$, then $xs=(sx)s$ implies $x=sx$ implies $s=e$,
which is a contradiction. Thus, $s\in C(L)$. If $(xy)s\ne x(ys)$ for
some $x,y\in L$, then $(xy)s=(x\cdot ys)s$ implies $xy=x\cdot ys$
implies $y=ys$ implies $s=e$, which is a contradiction. Thus $s\in
N(L)$. Therefore $s\in C(L)~,~s\in N(L)$ implies $s\in Z(L)$.
$\spadesuit$

\begin{myrem}
The result in Lemma~\ref{c:unic} is similar to that proved for
Moufang loops in \cite{phd42}.
\end{myrem}

\begin{mylem}\label{1:14}
In an LC-loop or RC-loop $(L,\cdot )$ with a unique non-identity
square $s$, the following are true.
\begin{enumerate}
\item $|s|=2$. \item $|x|=4$ or $|x|=2$. \item $s\in N_\lambda$ or
$s\in N_\rho$ and $s\in N_\mu$.
\end{enumerate}
\end{mylem}
{\bf Proof}\\ By hypothesis, $x^2=s$ for all $x\in L$.
\begin{enumerate}
\item Since $s\in L$, $s^2=s$ implies $s^{-1}s^2=s^{-1}s$ or
$s^2s^{-1}=ss^{-1}$ implies $s=e$. This is a contradiction, thus
$s^2=e$ if and only if $|s|=2$. \item $x^2=s$ implies
$x^4=x^2x^2=s^2=e$ implies $x^4=e$ or $x^2=e$.
\item In an LC-loop, $x^2\in N_\lambda ,N_\mu$, thus $s\in N_\lambda
,N_\mu$. In an RC-loop, $x^2\in N_\rho ,N_\mu$, thus $s\in N_\rho
,N_\mu$.
\end{enumerate} $\spadesuit$

\begin{mylem}\label{1:15}
Let $(L,\cdot)$ be a loop.
\begin{enumerate}
\item If $L$ is an LC-loop, $L$ has a unique non-identity square
$s$ if and only if $J=R_s^{-1}=R_{s^{-1}}^{-1}$ or $J=I$.
\item If $L$ is an RC-loop, $L$ has a unique non-identity square $s$
if and only if $J=L_s^{-1}=L_{s^{-1}}^{-1}$ or $J=I$.
\end{enumerate}
\end{mylem}
{\bf Proof}
\begin{enumerate}
\item $x^2=s$ if and only if $x^{-1}x^2=x^{-1}s$ if and only if
$x=x^{-1}s=xJR_s$ if and only if $JR_s=I$ if and only if
$J=R_s=R_{s^{-1}}^{-1}$.

$x^2=e$ if and only if $x=x^{-1}$ if and only if $x=xJ$ if and only
if $J=I$. \item $x^2=s$ if and only if $x^2x^{-1}=sx^{-1}$ if and
only if $x=sx^{-1}$ if and only if $x=xJL_s$ if and only if $I=JL_s$
if and only if $J=L_s^{-1}=L_{s^{-1}}^{-1}$.

$x^2=e$ if and only if $x=x^{-1}$ if and only if $x=xJ$ if and only
if $J=I$.
\end{enumerate} $\spadesuit$

\begin{myth}\label{1:16}
Let $(L,\cdot )$ be a L. I. P. (R. I. P. ) RC(LC)-loop with a unique
non-identity square $s$. The following are true.
\begin{enumerate}
\item $s\in Z(L,\cdot )$, hence $L$ is centrum square. \item
$J=L_s$($J=R_s$). \item $x^2y^2\ne (xy)^2\ne y^2x^2$ : $x\mapsto
x^2$ is neither an automorphism nor an anti-automorphism. \item
$(a,b,c)=(bc\cdot a)(ab\cdot c)$ :
\begin{description}
\item[(a)] $ab=a^{-1}b^{-1}$ if and only if $(J,J,I)\in AUT(L,\cdot
)$. \item[(b)] $(a,b,a)=(bs)(ab\cdot a)$ or $(a,b,a)=b(ab\cdot
a)$.
\end{description}
\item $L$ is a group or Steiner loop. \item If $L$ is not
commutative, $s$ is the unique non-identity commutator of $L$ for
the case of been a C-loop.
\end{enumerate}
\end{myth}
{\bf Proof}\\
\begin{enumerate}
\item $x^2=s$ implies $x=sx^{-1}$ implies $x^{-1}=s^{-1}x$.
Using the result in Lemma~2.1 of \cite{phd12},
$x^{-1}=(sx^{-1})^{-1}=(x^{-1})^{-1}s^{-1}$ implies
$x^{-1}=xs^{-1}$. Thus, $x^{-1}=s^{-1}x=xs^{-1}$ implies $sx=xs$
implies $s\in Z(L,\cdot )$. \item This follows from
Lemma~\ref{1:15}.
\item If $(xy)^2=x^2y^2$ or $(xy)^2=y^2x^2$ then $s=s^2$ implies
$s=e$ which is a contradiction. \item $(a,b,c)=[a(bc)]^{-1}\cdot
(ab)c=(bc)^{-1}a^{-1}\cdot (ab)c=(c^{-1}b^{-1})a^{-1}\cdot (ab\cdot
c)=[s^{-1}(bc)](s^{-1}a)\cdot (ab\cdot c)=(bc\cdot
s^{-1})(s^{-1}a)\cdot (ab\cdot c)=(bcs^{-2}\cdot a)(ab\cdot
c)=(bc\cdot a)(ab\cdot c)$. So, $(a,b,c)=(bc\cdot a)(ab\cdot c)$.
\begin{description}
\item[(a)] Put $c=e$, $(a,b,e)=(ba)(ab)=e$ if and only if
$(ab)=(ba)^{-1}$ if and only if $ab=a^{-1}b^{-1}$ if and only if
$(J,J,I)\in AUT(L,\cdot )$. \item[(b)] Let $c=a$, $(a,b,a)=(ba\cdot
a)(ab\cdot a)=(ba^2)(ab\cdot a)=(bs)(ab\cdot a)$ or $b(ab\cdot a)$.
\end{description}
\item This follows from Lemma~\ref{1:14}. \item
$(x,y)=x^{-1}y^{-1}\cdot xy=(x^{-1}y^{-1})(xy^{-1}\cdot
y^2)=((x^{-1}y^{-1})(xy^{-1})\cdot y^2=[x^{-2}(xy^{-1})\cdot
(xy^{-1})]y^2=x^{-2}[(xy^{-1})(xy^{-1})]y^2=e$ or $s$. Thus, $L$
is either commutative or else $s$ is its unique non-identity
commutator.

$(x,s)=x^{-1}s^{-1}\cdot xs=s$ implies $x^{-1}s\cdot xs=s$ implies
$x^{-1}R_s\cdot xR_s=s$ implies $xJ^2\cdot x^{-1}=s$ implies
$xx^{-1}=s$ implies $s=e$. This is a contradiction, so $(x,s)=e$
implies $s\in C(L,\cdot )$.
\end{enumerate} $\spadesuit$

\begin{mycor}\label{1:17}
A C-loop with unique non-trivial square is a group.
\end{mycor}
{\bf Proof}\\
By Lemma~\ref{1:14} and Theorem~\ref{1:16}, $L$ is central square
and is of exponent 4. Thus by Theorem~\ref{1:10}, $L$ is a group.
$\spadesuit$

\begin{myrem}
A C-loop with a unique non-trivial square is an A-loop.
\end{myrem}

\begin{myth}\label{1:18}
Let $(G,\cdot )$ and $(H,\circ )$ be two distinct loop such that
the triple $\alpha =(A,B,C)$ is an isotopism of $G$ upon $H$.
\begin{enumerate}
\item If $G$ is a central square C-loop of exponent $4$, then $H$ is
a C-loop and an A-loop. \item If $G$ is a C-loop with a unique
non-identity square, then $H$ is a C-loop and an A-loop.
\end{enumerate}
\end{myth}
{\bf Proof}\\
\begin{enumerate}
\item By Theorem~\ref{1:10}, $G$ is a group and since groups are
G-loops, $H$ is a group, thus both a C-loop and an A-loop. \item
The same argument follows by Corollary~\ref{1:17}.
\end{enumerate} $\spadesuit$

\begin{myrem}
In \cite{phd60} we considered only isotopes of central loops under
triples of the type $(A,B,B)$ and $(A,B,A)$. But
Theorem~\ref{1:18} addresses the general triple $(A,B,C)$ where
$A,B$ and $C$ are all distinct.
\end{myrem}

\begin{mycor}\label{1:19}
Let $(G,\cdot )$ and $(H,\circ )$ be distinct loops. If the triple
$(A,B,C)$ is an isotopism of $G$ upon $H$ such that
$(I,L_z^2,JL_z^2J),(R_z^2,I,JR_z^2J)\in AUT(G,\cdot )$ for all $z\in
G$, then $H$ is a C-loop and an A-loop.
\end{mycor}
{\bf Proof}\\
This follows by Theorem~\ref{c:centralsquare} and
Theorem~\ref{1:18}. $\spadesuit$

\begin{myth}\label{1:20}
Under  a triples of the form $(A,A,C)$, 'unique non-identity
square' is an isotopic invariant property for loops.
\end{myth}
{\bf Proof}\\
Let $(A,A,C)~:~(G,\cdot )\to (H,\circ )$ where $G$ and $H$ are two
distinct loops be an isotopism if and only if $xA\circ yA=(x\cdot
y)C$. Let $y=x$, $xA\circ xA=(xA)^2=(x\cdot x)C=x^2C$. If $s$ is the
unique non-identity square in $G$, i.e $x^2=s$ or $x^2=e$ for all
$x\in G$ then $s'=sC=(xA)^2=y'^2$ or $y'^2=(xA)^2=x^2C=eC=e'$ for
all $y'\in H$ with $e'$ as the identity element in $H$. So, $s'$ is
the unique non-identity square element in $H$. $\spadesuit$

\begin{mycor}\label{1:21}
Central loops with unique non-identity square are isotopic
invariant.
\end{mycor}
{\bf Proof}\\
This follows from Theorem~\ref{1:18} and Theorem~\ref{1:20}.
$\spadesuit$

\section{Weak and Cross inverse properties in Central Loops}
According to \cite{phd31}, the W. I. P. is a generalization of the
C. I. P.. The latter was introduced and studied by R. Artzy
\cite{phd45} and \cite{phd30} while the formal was introduced by J.
M. Osborn \cite{phd43} who also investigated the isotopy invariance
of the W. I. P.. Huthnance Jr. \cite{phd44}, proved that the
holomorph of a W. I. P. L. is a W. I. P. L.. A loop property is
called universal(or at times, a loop is said to be universal
relative to a particular property) if the loop has the property and
every loop isotope of such a loop possesses such a property. A
universal W. I. P. L. is called an Osborn loop. Huthnance Jr.
\cite{phd44} investigated the structure of some holomorph of Osborn
loops while Basarab \cite{phd46} studied Osborn loops that are
G-loops. Moufang and conjugacy closed loops have been found to be
Osborn loops in \cite{phd33} but in this section we shall
investigate this for C-loops.
\begin{myth}\label{1:22}
Let $(L,\cdot)$ be an LC(RC)-loop of exponent 3. $L$ is centrum
square if and only if $L$ is a C. I. P. L..
\end{myth}
{\bf Proof}\\
$x^3=e$ if and only if $x^2=x^{-1}$.

$x^2y=yx^2$ if and only if $x^{-1}y=yx^{-1}$ if and only if
$x(x^{-1}y)=x(yx^{-1})$ if and only if $y=x(yx^{-1})$ if and only if
the C. I. P. holds in the LC-loop $L$.

$x^2y=yx^2$ if and only if $yx^{-1}=x^{-1}y$ if and only if
$(yx^{-1})x=(x^{-1}y)x$ if and only if $y=(x^{-1}y)x$ if and only if
the C. I. P. holds in the RC-loop $L$. $\spadesuit$

\begin{myrem}
In the proof above, the fact in \cite{phd30} that in a C. I. P. L.,
the identity $xy\cdot xJ=y$ and $x(y\cdot xJ)=y$ are equivalent is
used. In fact, in a R. I. P. loop $L$, $L$ has the C. I. P. if and
only if $(xJ\cdot y)x=y$.
\end{myrem}

\begin{mycor}\label{1:23}
Let $(L,\cdot)$ be an LC(RC)-loop of exponent $3$. If $L$ is centrum
square, then
\begin{enumerate}
\item $L$ has both the A. I. P. and A. A. I. P. \item $L$ has the
W. I. P. \item $N=N_\lambda =N_\rho =N_\mu$. \item $n\in N$ implies
$n\in Z(L)$. \item $L$ is a commutative group.
\end{enumerate}
\end{mycor}
{\bf Proof}
\begin{enumerate}
\item By Theorem~\ref{1:22}, $L$ is a C. I. P. L. According to
\cite{phd30} and \cite{phd31}, a C. I. P. L. has the A. I. P.
Thus, the first part is true. The second part follows from the
fact that $x^2=x^{-1}$. \item This follows from the statement in
\cite{phd3} that W. I. P. is a generalization of C. I. P..
\end{enumerate}

(3) and (4) follow from \cite{phd31} and \cite{phd30} hence (5) is
as a result of the hypothesis. $\spadesuit$

\begin{myrem}
In Corollary~\ref{1:23}, $n\in Z(L)$ for all $n\in N$. Recall that
by hypothesis, $n=x^2\in N$ for all $x\in L$. Thus we expect $x^2\in
Z(L)$ which is actually the hypothesis $x^2\in C(L)$.
\end{myrem}

\begin{mycor}\label{1:23.5}
Let $(L,\cdot)$ be a C-loop of exponent $3$. If $L$ is central
square, then
\begin{enumerate}
\item $L$ has the W. I. P. and C. I. P. \item $L$ is a commutative
group.
\end{enumerate}
\end{mycor}
{\bf Proof}\\
A loop $(L,\cdot)$ is a C-loop if and only if it both an LC-loop and
an RC-loop. Thence, (1) follows from Theorem~\ref{1:22} and
Corollary~\ref{1:23}, while (2) follows from the latter.
$\spadesuit$

\begin{mylem}\label{1:24} LC(RC,C)-loop of exponent $3$ is a group.
\end{mylem}
{\bf Proof}\\
The proof of this is simple. $\spadesuit$

\paragraph{}
Despite the fact that central loops of exponent $3$ are groups it
will be interesting to study non-commutative central loops of
exponent 3 that have the C. I. P. since there exist groups that do
not have the C. I. P. From Theorem~\ref{1:22}, it can be seen that
the study of LC(RC)-loops of exponent $3$ with C. I. P. is
equivalent to the study of centrum square LC(RC)-loops of exponent
$3$.

The existence of central loops of exponent $3$ can be deduced from
\cite{phd56}, \cite{phd58} and \cite{phd10}. According to
\cite{phd58} and \cite{phd10}, the order of every element in a
finite LC(RC)-loop divides the order of the loop. Since $|x|=3$ for
all $x\in L$, $3\vert |L|$ which is possible because:

\begin{displaymath}
|L|=\left\{\begin{array}{lll} 2m, & \textrm{where}~m\ge3, &
\textrm{$L$
is non-left(right) Bol LC(RC)-loop}\\
4k, & \textrm{where}~k>2, & \textrm{$L$ is non-Moufang but both
lef t(right)-Bol and LC(RC)-loop}
\end{array}\right.
\end{displaymath}

The possible orders of finite RC-loops were proved in
\cite{phd10}. In \cite{phd10}, $m>3$ but in \cite{phd56}, the
author gave an example when $m=1$.

This type of central loops are flexible. This follows from
subsection~7.1. RC(LC)-loops are square flexible and $x\mapsto
x^{-1}$ is a permutation if and only if $x\mapsto x^2$ is a
permutation for a loop of exponent $3$. Also, if such finite central
loops are extra loops, they will be groups. This is as a result of
the claim in \cite{phd56} that an extra loop in which $x\mapsto x^2$
is a bijection is a group.

Flexible C-loops fall in a class of loops called ARIF loops and has
also been studied in \cite{phd22} and \cite{phd59}. Authors in
\cite{phd22} constructed a flexible C-loop that is not a RIF loop of
order $24$(this is divisible by $3$).

\subsection{Osborn Central-loops}

\begin{myth}\label{1:25}
Let $L$ be an LC(RC)-loop. $L$ has the R. I. P. (L. I. P.) if and
only if has the W. I. P.
\end{myth}
{\bf Proof} Let $(L,\cdot)$ be a LC-loop with the W. I. P.. Then for
all $x,y\in L,~y(xy)^\rho=x^\rho$. Let $xy=z$, then $x^\lambda
(xy)=x^\lambda z$ implies $y=x^\lambda z$, thus $(x^\lambda z)z^\rho
=x^\rho$ implies $(x^{-1}z)z^\rho =x^{-1}$. Replacing $x^{-1}$ with
$x$, $(xz)z^\rho =x$. Thence, $L$ has the R. I. P.

Conversely, if $L$ has the I. P., then $y(xy)^\rho
=y(xy)^{-1}=y(y^{-1}x^{-1})=x^{-1}=x^\rho$ hence it has the W. I. P.
Let $L$ be a RC-loop with the W. I. P. Then for all $x,y\in
L,~y(xy)^\rho =x^\rho$ if and only if $(xy)^\lambda\cdot
x=y^\lambda$. Let $xy=z$, then $(xy)y^\rho =zy^\rho$ implies
$x=zy^\rho$. Thus, $z^\lambda (zy^\rho )=y^\lambda$ implies
$z^\lambda (zy^{-1})=y^{-1}$. Replacing $y^{-1}$ by $y$, $z^\lambda
(zy)=y$. Thus, $L$ has the L. I. P..

Conversely, if $L$ has the I. P., then $y(xy)^\rho
=y(xy)^{-1}=y(y^{-1}x^{-1})=x^{-1}=x^\rho$, hence it has the W. I.
P. $\spadesuit$

\begin{mycor}\label{1:26}
Let $(L,\cdot)$ be an LC(RC)-loop with R. I. P.(L. I. P.). Then
\begin{enumerate}
\item $N(L)=N_\lambda(L)=N_\rho(L)=N_\mu(L)$. \item
$(I,R_{x^2},R_{x^2})\Big((L_{x^2},I,L_{x^2})\Big)\in AUT(L)$.
\item
$(L_x^2,R_{x^2},R_{x^2}L_x^2)\Big((L_{x^2},R_x^2,L_{x^2}R_x^2)\Big)\in
AUT(L)$.
\end{enumerate}
\end{mycor}
{\bf Proof}\\
By Theorem~\ref{1:25}, $L$ has the W. I. P.. According to
\cite{phd43}, in a W. I. P. L., the three nuclei coincide, so (1)
is true. Thus for an LC-loop, $x^2\in N_\rho$ and for an RC-loop,
$x^2\in N_\lambda$.

Hence for an LC-loop $L$, $(L_x^2,I,L_x^2),(I,R_{x^2},R_{x^2})\in
AUT(L)$ implies
$(L_x^2,R_{x^2},L_x^2R_{x^2})=(L_x^2,R_{x^2},R_{x^2}L_x^2)\in
AUT(L)$.

For an RC-loop $L$, $(I,R_x^2,R_x^2),(L_{x^2},I,L_{x^2})\in AUT(L)$
implies
$(L_{x^2},R_x^2,R_x^2L_{x^2})=(L_{x^2},R_x^2,L_{x^2}R_x^2)\in
AUT(L)$.

So, (2) and (3) are true. $\spadesuit$

\begin{myrem}
Corollary~\ref{1:26} is true for left(right) Bol loops(i.e
LB(RB)-loops). This can be seen from the fact in \cite{phd3} that
a RB(LB)-loop has the L. I. P.(R. I. P.) if and only if it is a
Moufang loop which is obviously a W. I. P. L. as mentioned in
\cite{phd33}.
\end{myrem}

The following facts which is true in Bol loops is also true in
central loops.

\begin{myth}\label{1:27}
A LC(RC)-loop $(L,\cdot)$ is a C-loop if and only if any of the
following equivalent statements about $L$ is true.
\end{myth}
\begin{enumerate}
\item $L$ has the R. I. P.(L. I. P.) \item $L$ has the R. A. P.(L.
A. P.) \item $L$ is a RC(LC)-loop. \item $L$ has the A. A. I. P..(
i.e $(xy)^{-1}=y^{-1}x^{-1}$) \item $L$ has the W. I. P..
\end{enumerate}
{\bf Proof}\\  A C-loop satisfies (1)-- (2). The proof of the
converse is as follows.
\begin{enumerate}
\item $L$ is an LC-loop if and only if $(x\cdot xy)z=x(x\cdot
yz)$ if and only if $[(x\cdot xy)z]^{-1}=[x(x\cdot yz)]^{-1}$ if and
only if $z^{-1}(x\cdot xy)^{-1}=(x\cdot yz)^{-1}x^{-1}$ if and only
if $z^{-1}((xy)^{-1}\cdot x^{-1})=((yz)^{-1}\cdot x^{-1})x^{-1}$ if
and only if $z^{-1}(y^{-1}x^{-1}\cdot x^{-1})=(z^{-1}y^{-1}\cdot
x^{-1})x^{-1}$. Carrying out the replacement $x,y,z$ for
$x^{-1},y^{-1},z^{-1}$ we have $z(yx\cdot x)=(zy\cdot x)x$ if and
only if $L$ is an RC-loop. Hence, $L$ is a C-loop.
\item According to \cite{phd58}, $L$ is an LC-loop if and only if
$x\cdot(y\cdot yz)=(x\cdot yy)z$ for all $x,y,z\in L$ while $L$ is
an RC-loop if and only if $(zy\cdot y)x=z(yy\cdot x)$ for all
$x,y,z\in L$. $x\cdot(y\cdot yz)=(x\cdot yy)z$ if and only if
$x\cdot zL_y^2=xR_{y^2}\cdot z$ if and only if
$(R_{y^2},L_y^{-2},I)\in AUT(L,\cdot )$ for all $y\in L$ while
$(zy\cdot y)x=z(yy\cdot x)$ if and only if $zR^2\cdot x=z\cdot
xL_{y^2}$ if and only if $(R_y^2,L_{y^2}^{-1},I)\in AUT(L,\cdot)$
for all $y\in L$.

If $L$ has the right(left) alternative property,
$(R_y^2,L_y^{-2},I)\in AUT(L,\cdot )$ for all $y\in L$ if and only
if $L$ is a C-loop.
\item This is shown
in \cite{phd56}. \item This is equivalent to (1) : If $L$ has the L.
I. P.(R. I. P.) then $L$ has the R. I. P.(L. I. P.) implies the A.
A. I. P.. Conversely, if L. I. P. holds, then with $z=xy$, we have
$y=x^{-1}z$ which gives $z^{-1}=(x^{-1}z)^{-1}x^{-1}$ implies
$z^{-1}=(x^{-1}z)^{-1}x^{-1}$ implies $z^{-1}=(z^{-1}x)x^{-1}$.
Replacing $z^{-1}$ by $z$, we have $z=(zx)x^{-1}$ implies the R. I.
P.

Similarly, if $L$ has the R. I. P.(L. I. P.) then $L$ has the L. I.
P.(R. I. P.) implies the A. A. I. P.. Conversely, if R. I. P. holds,
then with $z=xy$, we have $x=zy^{-1}$. Thus,
$z^{-1}=y^{-1}(zy^{-1})^{-1}=y^{-1}(yz^{-1})$. Replacing $z^{-1}$ by
$z$, $z=y^{-1}(yz)$ implies L. I. P. \item This follows from (1) and
Theorem~\ref{1:25}.
\end{enumerate} $\spadesuit$

\begin{myth}\label{1:28}(\cite{phd33})
The following equivalent conditions define an Osborn loop
$(L,\cdot)$.
\begin{enumerate}
\item $x(yz\cdot x)=(x\cdot yE_x)\cdot zx$. \item $(x\cdot
yz)x=xy\cdot (zE_x^{-1}\cdot x)$. \item $(A_x,R_x,R_xL_x)\in
AUT(L)$. \item $(L_x,B_x,L_xR_x)\in AUT(L)$.
\end{enumerate}
where $A_x=E_xL_x,B_x=E_x^{-1}R_x$ and
$E_x=R_xL_xR_x^{-1}L_x^{-1}$.
\end{myth}

\begin{myth}\label{1:29}
Let $(L,\cdot)$ be a RC(LC)-loop. If $L$ has the L. I. P.(R. I. P.),
then $L$ then $L$ is an Osborn loop if  every element in $L$ is a
square.
\end{myth}
{\bf Proof}\\
\begin{enumerate}
\item Let $L$ be an RC-loop with L. I. P. then by
Theorem~\ref{1:25}, $L$ has the W. I. P..

$(A_{x^2},I,L_{x^2})\in AUT(L)$ if and only if $yA_{x^2}\cdot
z=(yz)L_{x^2}$ if and only if $yE_{x^2}L_{x^2}\cdot z=(yz)L_{x^2}$
if and only if $yR_{x^2}L_{x^2}R_{x^2}^{-1}L_{x^2}^{-1}L_{x^2}\cdot
z=(yz)L_{x^2}$ if and only if $yR_{x^2}L_{x^2}R_{x^2}^{-1}\cdot
z=(yz)L_{x^2}$ if and only if $yR_x^2L_{x^2}R_{x^2}^{-1}\cdot
z=(yz)L_{x^2}$ if and only if $yL_{x^2}R_x^2R_{x^2}^{-1}\cdot
z=(yz)L_{x^2}$ if and only if $yL_{x^2}\cdot z=(yz)L_{x^2}$ if and
only if $(L_{x^2}, I, L_{x^2})\in AUT(L)$ for all $x\in L$ which is
true by Corollary~\ref{1:26}.

Thus, $(I,R_x^2,
R_x^2)(A_{x^2},I,L_{x^2})=(A_{x^2},R_{x^2},R_{x^2}L_{x^2})\in
AUT(L)$. Using Theorem~\ref{1:28}, $L$ is an Osborn loop if every
element in $L$ is a square(i.e $y=x^2$ for all $y\in L$). \item Let
$L$ be an LC-loop. If $L$ has the R. I. P., then by Theorem
~\ref{1:25}, $L$ has the W. I. P..

$(I,B_{x^2},R_{x^2})\in AUT(L)$ if and only if $y\cdot
zB_{x^2}=(yz)R_{x^2}$ if and only if $y\cdot
zE_{x^2}^{-1}R_{x^2}=(yz)R_{x^2}$ if and only if $y\cdot
z(R_{x^2}L_{x^2}R_{x^2}^{-1}L_{x^2}^{-1})^{-1}R_{x^2}=(yz)R_{x^2}$
if and only if $y\cdot
zL_{x^2}R_{x^2}L_{x^2}^{-1}R_{x^2}^{-1}R_{x^2}=(yz)R_{x^2}$ if and
only if $y\cdot zL_{x^2}R_{x^2}L_{x^2}^{-1}=(yz)R_{x^2}$ if and only
if $y\cdot zR_{x^2}L_x^2L_{x^2}^{-1}=(yz)R_{x^2}$ if and only if
$y\cdot zR_{x^2}=(yz)R_{x^2}$ if and only if $(I,R_{x^2},R_{x^2})\in
AUT(L)$ for all $x\in L$ which is true by Corollary~\ref{1:26}.

Thus, $(L_x^2,I,L_x^2)(I,B_{x^2},R_{x^2})=(L_{x
^2},B_{x^2},L_{x^2}R_{x^2})\in AUT(L)$. Using Theorem~\ref{1:28},
$L$ is an Osborn loop if every element in $L$ is a square(i.e
$y=x^2$ for all $y\in L$).
\end{enumerate}
$\spadesuit$

\begin{mycor}\label{1:30}
Let $(L,\cdot)$ be a LC(RC)-loop with R. I. P.(L. I. P.). $L$ is an
Osborn loop if every element in $L$ is a square. Hence, $L$ is a
group.
\end{mycor}
{\bf Proof}\\
This follows from Theorem~\ref{1:29}. The last conclusion is as a
consequence of the fact that $x^2\in N(L)$.  $\spadesuit$

\begin{mycor}\label{1:31}
Let $(L,\cdot)$ be a C-loop. $L$ is an Osborn loop if every element
in $L$ is a square. Hence, $L$ is a group.
\end{mycor}
{\bf Proof}\\
This follows from Theorem~\ref{1:29}. The last conclusion is as a
consequence of the fact that $x^2\in N(L)$. $\spadesuit$

\begin{myrem}\label{1:32}
Corollary~\ref{1:31} informs us that a class of C-loops that are
Osborn are groups. Since the condition attached is not necessary and
sufficient, there is the possibility of the existence of C-loops
that are Osborn but non-associative.
\end{myrem}

\section{Question}
Does there exist a C-loop that is Osborn but non-associative, non
Moufang and non-conjugacy closed.

\section{Acknowledgement} The second author would
like to express his profound gratitude to the Swedish
International Development Cooperation Agency (SIDA) for the
support for this research under the framework of the Associateship
Scheme of the Abdus Salam International Centre for
theoretical Physics, Trieste, Italy. ~\\

~\\
~\\
~\\
\begin{enumerate}
\item Department of Mathematics,\\
Obafemi Awolowo University,\\ Il\'e If\`e, Nigeria.\\
{\bf e-mail}: jaiyeolatemitope@yahoo.com \item Department of Mathematics,\\
University of Ab\d e\'ok\`uta,\\ Ab\d e\'ok\`uta 110101, Nigeria.\\
{\bf e-mail}: ekenedilichineke@yahoo.com
\end{enumerate}
\end{document}